%% file: smo_pareto_hongahnchoi_arxiv.tex
\documentclass[twocolumn]{my_svjour3}         
\smartqed  

\usepackage{graphicx}
\usepackage{epsfig} 
\usepackage{epstopdf}
\usepackage{algorithm}
\usepackage{algorithmic}
\usepackage{nomencl}
\usepackage[cmex10]{amsmath, mathtools}
\usepackage{amsfonts, amssymb}
\usepackage{theorem}
\usepackage{setspace}
\usepackage{subfigure}
\usepackage{multirow}
\usepackage[usenames]{color}
\usepackage{array}
\usepackage[sort&compress, numbers]{natbib}
\usepackage{ifpdf}
\usepackage{balance}
\usepackage{enumitem}
\usepackage{times}

\newcommand{\be}{\begin{equation}}
\newcommand{\ee}{\end{equation}}
\newcommand{\beas}{\begin{eqnarray*} }
 \newcommand{\eeas}{\end{eqnarray*}}

\theorembodyfont{\rmfamily}
\theoremstyle{plain}\newtheorem{Rem}{Remark}
\theoremstyle{plain}
\theoremstyle{plain}
\theoremstyle{plain}
 \theoremstyle{break}

\theoremstyle{plain}
\newtheorem{Prop}{Proposition}
\theoremstyle{plain}

\journalname{Preprint submitted to \textit{Structural and Multidisciplinary Optimization}}
\begin{document}

\title{ Pruning-Based Pareto Front Generation for Mixed-Discrete Bi-Objective Optimization
}


\author{SeungBum Hong         \and
           Jaemyung Ahn \and
           Han-Lim Choi$^\ast$\thanks{$^\ast$Corresponding Author. Email:hanlimc@kaist.ac.kr; Tel:+82-42-350-3727; Fax:+82-42-350-3710}
}


\institute{S. Hong \and J. Ahn \and H.-L. Choi
\at
Division of Aerospace Engineering, KAIST \\
 291 Daehak-ro, Yuseong, Daejeon 305-701, Korea. \\
              }

\date{}

\maketitle

\begin{abstract}
This note proposes an effective pruning-based Pareto front generation method in mixed-discrete bi-objective optimization. The mixed-discrete problem is decomposed into multiple continuous subproblems; two-phase pruning steps identify and prune out non-contributory subproblems to the Pareto front construction. The efficacy of the proposed method is demonstrated on two benchmark examples.
\end{abstract}


\section{Introduction}

Consider a bi-objective optimization (BOO) problem whose design vector ($\mathbf{x}$) has both of continuous ($\mathbf{y} = [y_1,\cdots,y_{n_y}]$) and discrete ($\mathbf{z}=[z_1,\cdots,z_{n_z}]$) components:
\begin{align}
& \min_{\mathbf{x}} \mathbf{J}(\mathbf{x}) =
\min_{[\mathbf{y}~\mathbf{z}]} \mathbf{J} ([\mathbf{y}~\mathbf{z}]) =
\begin{bmatrix}   J_1 (\mathbf{y},\mathbf{z}) &  J_2 (\mathbf{y},\mathbf{z}) \end{bmatrix}^\top \label{MOO_master}  \tag{\textbf{P}}  \\
& \text{subject to} \notag \\
&\mathbf{g}(\mathbf{y},\mathbf{z}) \leq 0, \qquad \mathbf{h}(\mathbf{y},\mathbf{z}) = 0, \notag \\
&y_i \in [l_i, u_i], \qquad i=1,\cdots,n_y, \notag \\
&z_j \in Z_j = \left\{z_j^{1},\cdots,z_j^{|Z_j|} \right\},\qquad j=1,\cdots,n_z,\notag
\end{align}
where $\mathbf{g}$ is the inequality constraint vector, $\mathbf{h}$ is the equality constraint vector, $l_i$ and $u_i$ are the lower and upper bounds of the $i^{\text{th}}$ continuous design variable ($y_i$), and $Z_j$ is the set of values that $j^{\text{th}}$  discrete design variable ($z_j$) can take. Let $\mathcal{X}^\star$ be the set of design vectors that are Pareto optimal solutions of \ref{MOO_master}:
\begin{equation}
\mathcal{X}^\star=\left\{\mathbf{x}^\star \in \mathcal{X}| \nexists \mathbf{x} \in \mathcal{X} \setminus \{ \mathbf{x}^\star \}~\text{s.t.}~ \mathbf{J} (\mathbf{x}) \leq \mathbf{J} (\mathbf{x}^\star) \right\} \label{pareto_opt}
\end{equation}
where $\mathcal{X}$ is the set of feasible design vectors. The problem of \textit{Pareto front generation} is equivalent to determining $\mathcal{X}^\star$.

In case the objective function and the constraints are linear, resulting in a multi-objective mixed-integer program, several tailored algorithms have been proposed~\citep{Alves_EJOR00, Ulungu94}. For nonlinear mixed-discrete problems,  meta-heuristic approaches such as genetic algorithm~\citep{Deb_TEC00}, evolutionary programming~\citep{Meza_SMC09}, particle swarm optimization~\citep{Wang07}, and tabu search~\citep{Sendin_SB09} have often been adopted for Pareto front generation.  As deterministic approach, an iterative two-phase procedure that solves given number of mixed-integer nonlinear programs and then a sequence of continuous nonlinear programs (NLPs) in \citet{Mela07} is the only work reported in the literature to the authors' best knowledge. A common fact for all the previous work on nonlinear cases is that some number of mixed-integer nonlinear programs need to be solved; this requirement might be an issue for practical purposes (in particular when the discrete variables are categorical). On the contrary, this note takes advantage of decomposition and pruning methodology that does not require solution of complex mixed-discrete nonlinear programs.

\section{Subproblem Decomposition and Pruning}\label{method}
\subsection{Approach}
One way to generate the Pareto front of the original BOO is to divide \ref{MOO_master} into subproblems with specific discrete design vectors and construct $\mathcal{X}^\star$ by systematically synthesizing the solutions of the subproblems. First, define the set of discrete design vectors, $\mathcal{Z}\triangleq Z_1 \times \cdots \times Z_{n_z}$, and associated index set $\mathcal{K} = \{ 1, 2, \dots, |\mathcal{Z}| \}$. Let $\mathbf{z}_k,~k \in \mathcal{K} $ be the $k^{\text{th}}$ element of $\mathcal{Z}$; a subproblem of \ref{MOO_master} associated with this discrete realization, denoted as \ref{MOO_zk}, can be defined as:
\begin{align}
&\min_{\mathbf{y}} \mathbf{J} ([\mathbf{y}~\mathbf{z}_k]) =
\begin{bmatrix}   J_1 (\mathbf{y},\mathbf{z}_k) &  J_2 (\mathbf{y},\mathbf{z}_k) \end{bmatrix}^\top
\label{MOO_zk}  \tag{\textbf{P}$_k$}  \\
& \text{subject to} \notag \\
&\mathbf{g}(\mathbf{y},\mathbf{z}_k) \leq 0, \qquad \mathbf{h}(\mathbf{y},\mathbf{z}_k) = 0, \notag \\
&y_i \in [l_i, u_i], \qquad i=1,\cdots,n_y. \notag
\end{align}
The set of Pareto optimal solutions for \ref{MOO_zk} is defined as
\begin{equation*}
\mathcal{X}_k^\star =\left\{ \mathbf{x}^\star \in \mathcal{Y}_k \times \{\mathbf{z}_k \} |
 \nexists~ \mathbf{y} \in \mathcal{Y}_k ~\text{s.t.}~
 \mathbf{J} ([\mathbf{y}, \mathbf{z}_k] ) \leq \mathbf{J} (\mathbf{x}^\star) \right\}
\end{equation*}
where $\mathcal{Y}_k$ is the set of feasible continuous design vectors of the subproblem.

$\mathcal{X}_k^\star$ can be obtained relatively easily using normal boundary intersection (NBI)~\citep{DasDennis_JO98, MarlerArora_SMO04, Motta_SMO12} or the weighted sum (WS) method~\citep{Hwang1979, MarlerArora_SMO04, Marler_SMO10, Kim_SMO05} combined with reliable nonlinear programming (NLP) solvers. One brute-force way of obtaining $\mathcal{X}^\star$ is to first compute $\mathcal{X}_k^\star$ for all possible discrete realization $\mathbf{z}_k$ and then identify, among those subproblem solutions, design vectors satisfying (\ref{pareto_opt}). But this approach can be computationally intractable if the discrete design space is very large, i.e., large $|\mathcal{Z}|$.

Note that the Pareto optimal solutions for some subproblems may have no common elements with $\mathcal{X}^\star$, while the others have common elements with $\mathcal{X}^\star$ and thus contribute to constructing the Pareto front of \ref{MOO_master}.  Define the index set of irrelevant ($\mathcal{K}_\emptyset$) and relevant ($\mathcal{K}_1$) subproblems, respectively:
\begin{equation*}
\mathcal{K}_{\emptyset} = \left\{k\in \mathcal{K} |\mathcal{X}^\star_k \cap \mathcal{X}^\star = \emptyset \right\}, ~~
\mathcal{K}_{1} = \mathcal{K} \setminus \mathcal{K}_{\emptyset} \\
\end{equation*}
Then, the Pareto optimal solution to \ref{MOO_master} can be obtained by collecting non-dominated solutions out of the relevant Pareto subproblem solutions:
\begin{equation*}
\mathcal{X}^\dag=\left\{\mathbf{x}^\dag \in \mathcal{X}_{\mathcal{K}_{1}}^\star |  \nexists \mathbf{x} \in \mathcal{X}_{\mathcal{K}_{1}}^\star \setminus \{ \mathbf{x}^\dag \} ~\text{s.t.}~ \mathbf{J} (\mathbf{x}) \leq \mathbf{J} (\mathbf{x}^\dag)
\right\} \label{pareto_opt_prune}
\end{equation*}
where $  \mathcal{X}_{\mathcal{K}_{1}}^\star = \bigcup_{k \in \mathcal{K}_1} \mathcal{X}_k^\star $.

Therefore, if $\mathcal{K}_1$ (or equivalently $\mathcal{K}_{\emptyset}$) can be identified in advance by some efficient procedure, computational complexity of solving \ref{MOO_master} will be significantly reduced, in particular, when $|\mathcal{K}_1| \ll |\mathcal{Z}|$.  This work proposes a set of heuristics to approximately identify $\mathcal{K}_1$ by solving constant number of nonlinear programs for each \ref{MOO_zk}.

\subsection{Algorithm}  \label{sec:algorithm}
\begin{figure}[t]
\centerline{\includegraphics[trim=242 167 220 193,clip,width=.95\columnwidth]{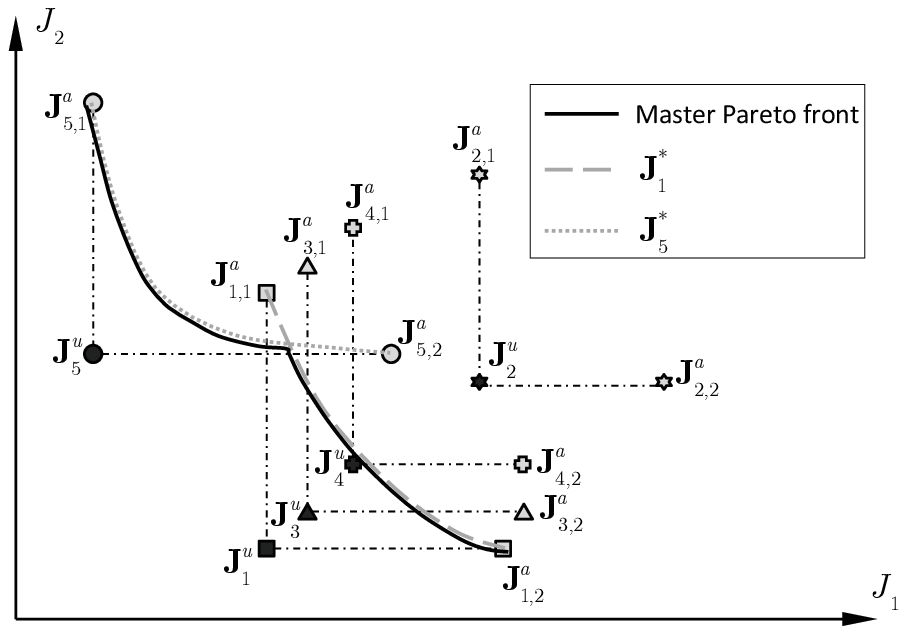}} \vspace*{-.0in} \caption{Illustration of Phase A Pruning}
  \label{fig:utopia}
  \vspace*{.1in}
\centerline{\includegraphics[trim=242 167 220 193,clip,width=.95\columnwidth]{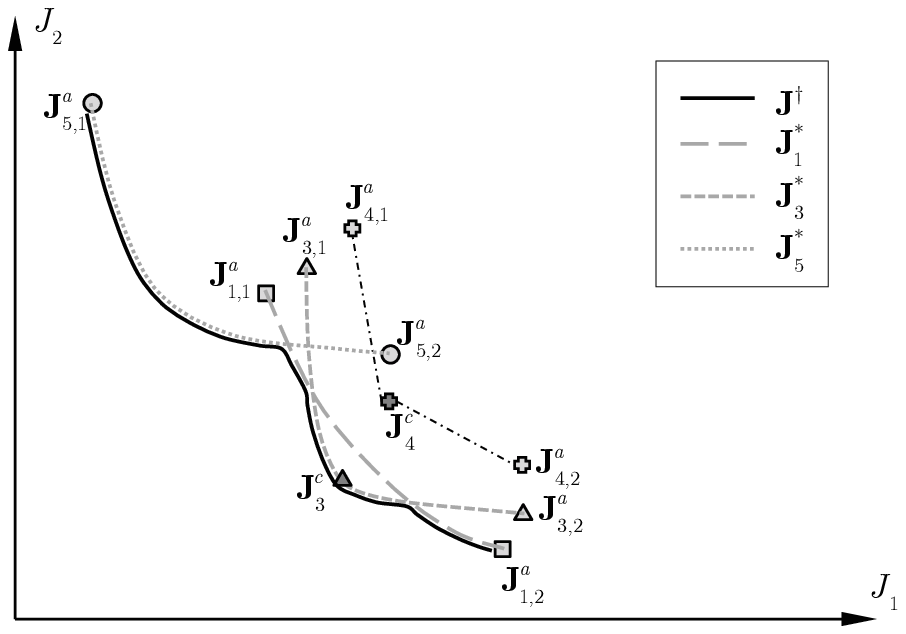}} \vspace*{-.0in} \caption{Illustration of Phase B Pruning}
  \label{fig:midpoint}
\vspace*{-.1in}
\end{figure}

This note presents a mechanism to prune a set of subproblems that are expected not to contribute to construction of the Pareto front of \ref{MOO_master}. The procedure consists of two phases: Phase A based on dominance of subproblem utopia points followed by Phase B based on dominance of center points of subproblem Pareto front.

\paragraph{Phase A-1: Computing subproblem anchor/utopia points}
The anchor points of \ref{MOO_zk} are obtained as
\begin{equation*}
\mathbf{J}_{k,i}^a= \mathbf{J}( \mathbf{y}_{k,i}^\sharp)  \triangleq [ J_{i,k}^{a,1}, J_{i,k}^{a,2} ]^\top ,  ~~i = \{ 1, 2\}
\end{equation*}
where $ \mathbf{y}_{k,i}^\sharp = \arg \min_{\mathbf{y} \in \mathcal{Y}_k} J_i (\mathbf{y}, \mathbf{z}_k)$, i.e., solution to the two sole-objective optimization problems. The utopia point of \ref{MOO_zk} is then be computed as
\begin{equation*}
\mathbf{J}_k^u = [J_{1,k}^{a,1}, J_{2,k}^{a,2}]^\top.
\end{equation*}
This step computes the anchor points and  the utopia points for all \ref{MOO_zk} (see five sets of anchor/utopia points Fig. \ref{fig:utopia}).

\paragraph{Phase A-2: Generating a master Pareto front}
Cross-checking of dominance between the utopia points allows for identification of $\mathcal{K}_{1}^m$ that can be used to compute an approximate Pareto front:
\begin{equation}
\mathcal{K}_1^m = \left\{k |\nexists~ l \in  \mathcal{K} \setminus \{k\},~ \mathbf{J}_l^u \leq \mathbf{J}_k^u \right\}. \label{eq:k1m}
\end{equation}
Once $\mathcal{K}_1^m$ is determined, a Pareto front with this subproblem set can be obtained, $\mathcal{X}_{\mathcal{K}_1^m}^\star$, which is termed \textit{master} front herein. For example, in Fig. \ref{fig:utopia}, utopia points for $k=1, 5$ are non-dominated; a master Pareto front is generated by obtaining the solutions of subproblems $\textbf{P}_1$ and $\textbf{P}_5$ and selecting non-dominated elements.

\paragraph{Phase A-3: Pruning irrelevant subproblems}
For the subproblems not considered in construction of the master front, dominance of those utopia points compared to the master front is investigated to obtain
\begin{equation*}
\mathcal{K}_{\emptyset}^u = \left\{ k \left| \exists \mathbf{x} \in \mathcal{X}^\star_{\mathcal{K}_1^m}~\text{s.t.}~
\mathbf{J}_k^u \geq \mathbf{J}(\mathbf{x}) \right. \right\}, \label{eq:k0u}
\end{equation*}
where $\mathcal{X}^\star_{\mathcal{K}_1^m} \triangleq \bigcup_{k \in \mathcal{K}_1^m} \mathcal{X}_k^\star$; the subproblems in this set $\mathcal{K}_{\emptyset}^u$ are pruned.
(See in Fig. \ref{fig:utopia} the utopia point for $\textbf{P}_2$ is dominated by the master front.)

As a result, at the end of Phase A, subproblems in the set
\begin{equation}
\mathcal{K}_1^u = \mathcal{K} \setminus \mathcal{K}_{\emptyset}^u  \label{eq:k1u}
\end{equation}
are left for consideration of Pareto front generation of \ref{MOO_master}, some of whose subproblem Pareto fronts have already been created in Phase A-2.

\begin{Prop} \label{prop:phase_a}
The Pareto front constructed with the subproblem set $\mathcal{K}_1^u$ is identical to the true Pareto front $\mathcal{X}^\star$.
\begin{proof}
It suffices to prove that any $\mathbf{x}_l \in \mathcal{X}_l^\star,~~l \in \mathcal{K}_{\emptyset}^u$ is dominated by some other design vector; thus, $\mathcal{X}_l^\star \cap \mathcal{X}^\star = \emptyset$.  For such $\mathbf{x}_l$,
$
\mathbf{J}(\mathbf{x}_l) \geq \mathbf{J}_l^u
$
by the definition of utopia point. The fact that $\mathbf{J}_l^u$ is pruned implies that
$
\exists \mathbf{x}_k \in \mathcal{X}_{\mathcal{K}_1^m}~~\text{s.t.}~~   \mathbf{J}_l^u \geq \mathbf{J}(\mathbf{x}_k)
 $
Therefore, $\mathbf{x}_l$ is dominated by at least one element in the master front; thus, it cannot be included in $\mathcal{X}^\star$. \qed
\end{proof}
\end{Prop}

Depending on the problem type, $\mathcal{K}_1^u$ may not be substantially smaller than $\mathcal{K}$; in this case, the following Phase B procedure can improve the computational efficiency.

\paragraph{Phase B-1: Computing approximate subproblem center points}
For subproblems that have neither been pruned through Phase A nor used to build the master Pareto front, identify one point on the subproblem Pareto front. This step solves one nonlinear program per subproblem, in particular if the weighted sum method is adopted for Pareto front calculation, the following NLP is solved:
\begin{equation*}
\mathbf{y}_k^c = \arg \min_{\mathbf{y}} 0.5  J_1 (\mathbf{y}, \mathbf{z}_k) + 0.5 J_2 (\mathbf{y}, \mathbf{z}_k)
\end{equation*}
subject to the constraints for \ref{MOO_zk}; $\mathbf{x}_k^c = [\mathbf{y}_k^c~ \mathbf{z}_k]$ and $\mathbf{J}_k^c = \mathbf{J}(\mathbf{x}_k^c)$ are also computed accordingly.

\paragraph{Phase B-2: Pruning dominated center point}
A Pareto front for a subproblem can be approximated as an arc passing through the two anchor points and the center point; this step assess the dominance of the subproblem Pareto front based on the piecewise linear segment consisting of these three points. This step checks whether or not the center point is dominated by the master front constructed in Phase A-2; if dominated the associated subproblem is pruned from the candidate list of relevant subproblems:
\begin{equation*}
\mathcal{K}_{\emptyset}^c = \left\{ k \left| \exists \mathbf{x} \in \mathcal{X}^\star_{\mathcal{K}_1^m}~\text{s.t.}~
\mathbf{J}_k^c \geq \mathbf{J}(\mathbf{x}) \right. \right\}. \label{eq:k0c}
\end{equation*}
Observe in Fig. \ref{fig:midpoint} that $\textbf{P}_4 $ is pruned as $\mathbf{J}_4^c$ is dominated by the master front.

\paragraph{Phase B-3: Generating Pareto front with remaining subproblems} Then, the remaining index set becomes
\begin{equation}
\mathcal{K}_1^c = \mathcal{K}_1^u \setminus \mathcal{K}_{\emptyset}^c, \label{eq:k1c}
\end{equation}
and typical Pareto generation techniques such as weighted sum or NBI method can be used to construct the Pareto front of \ref{MOO_master}. As a result, the Pareto front is calculated using the subproblem Pareto fronts in $\mathcal{K}_1^c$. In Fig. \ref{fig:midpoint}, $\textbf{P}_3$ is included to construct the Pareto front. It should be pointed out that since $\mathcal{K}_1^c$ is not a superset of $\mathcal{K}_1$, optimality of the resulting Pareto front is not guaranteed with Phase B, in contrast to preservation of optimality only with Phase A.

\begin{Rem}
Throughout this two-phase process, the total number of NLPs to be solved to construct the Pareto front is
\begin{equation*}
N_{\text{AB}} = \underbrace{2 |\mathcal{K}|}_{\text{A-1}} + \underbrace{\beta |\mathcal{K}_1^m |}_{\text{A-2}} +  \underbrace{( |\mathcal{K}_1^u | -  |\mathcal{K}_1^m| )}_{\text{B-1}} + \underbrace{\beta (|\mathcal{K}_1^c| - |\mathcal{K}_1^m| )}_{\text{B-3}}
\end{equation*}
where $\beta$ is the number of points in each subproblem Pareto (the corresponding phase numbers are given beneath the underbraces). Note that without pruning $N_{\text{O}}= \beta |\mathcal{K}|$ nonlinear programs need to be solved; thus,
\begin{equation*}
\frac{N_{\text{AB}}}{N_{\text{O}}}  \leq  \frac{2}{\beta} + \frac{ |\mathcal{K}_1^u|}{  \beta |\mathcal{K}|}  + \frac{|\mathcal{K}_1^c|}{|\mathcal{K}|} \approx
 \begin{cases}
 \frac{3}{\beta} +  1 & |\mathcal{K}_1^c| \lesssim |\mathcal{K}| \\
 \frac{3}{\beta} & |\mathcal{K}_1^c| \ll |\mathcal{K}_1^u| \lesssim |\mathcal{K}| \\
 \frac{2}{\beta} & |\mathcal{K}_1^u| \ll |\mathcal{K}|
 \end{cases}
\end{equation*}
 Typically $\beta \sim \mathcal{O}(10)$ to ensure sufficient accuracy of the Pareto front; the pruning method can achieve $\mathcal{O}(10)$ times efficiency for typical cases at the price of $\mathcal{O}(0.1)$ overhead computation time in the worst case.
 \qed
\end{Rem}


\section{Numerical Examples}  \label{sec:numerics}
Two numerical examples are considered to demonstrate the effectiveness of the proposed algorithm.

\paragraph{Van Veldhuizen's Test Problem} One of the Van Veldhuizen's test suite that is known to exhibit non-trivial Pareto optimal set is considered \citep{Huband06}. The formulation is as the following:
\begin{align*}
& \min\left[\begin{array}{c}
J_{1}(\mathbf{x})\\
J_{2}(\mathbf{x})
\end{array}\right]=\left[\begin{array}{c}
\sum_{i=1}^{2}-10e^{-0.2\sqrt{x_{i}^{2}+x_{i+1}^{2}}}\\
\sum_{i=1}^{3}\left\{ |x_{i}|^{0.8}+5\sin(x_{i}^{3})\right\}
\end{array}\right] \label{MOP4}  \tag{\textbf{E}$_1$} \\
& \text{subject to} \\
& x_1 \in [-5,5], \qquad  x_2, x_3  \in \{-5,-4, \dots, 4, 5 \}. 
\end{align*}
 While $x_1$ is continuous, $x_2$, $x_3$ are discrete variables each of which can take 11 possible discrete values; $|\mathcal{K}| = 121$. For comparison, exhaustive search is implemented to obtain the true Pareto optimal points using weighted sum method. Table \ref{tab:res_summary} indicates that the proposed method identifies the same set of discrete realizations as the the true one; significant number of subproblems are pruned out by Phase A.

\begin{figure}[t]
\vspace*{-.05in}
  \centering
  \includegraphics[trim=5 5 5 5,clip,width=.7\columnwidth]{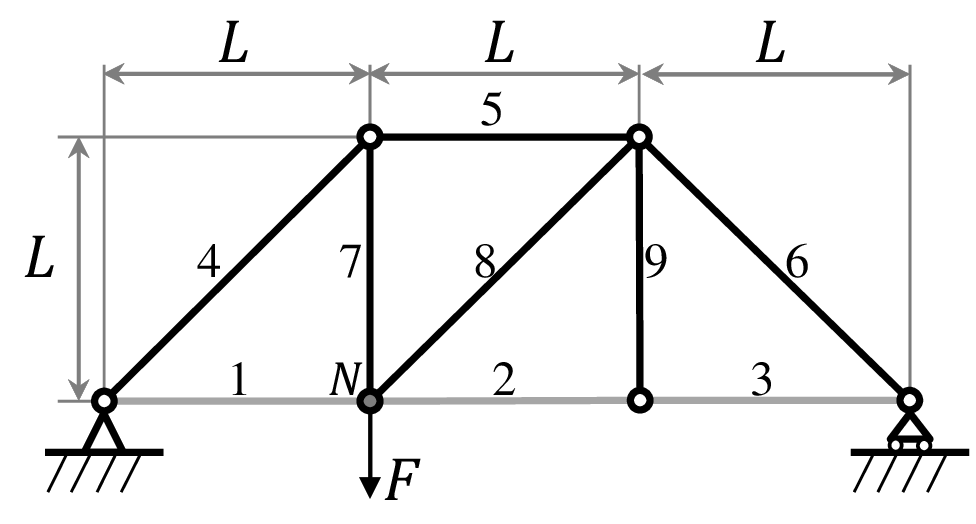}
  \caption{Illustration of the Nine bar truss.\label{fig:nineBar}}
\end{figure}

\begin{table}[t]
  \caption{Effects of pruning on number of subproblems} \vspace*{-.05in}
  \begin{center}
  \begin{tabular}{c c c c c}
   \hline
Case &  $|\mathcal{K}|$ & $|\mathcal{K}_1|$ & $|\mathcal{K}_1^u |$ & $|\mathcal{K}_1^c | $ \\
    \hline\hline
$   \textbf{E}_1$ & 121 &  4 & 5 & 4 \\
  $ \textbf{E}_2$ & 4,096 & 72 &  907 & 72 \\
     \hline
  \end{tabular}
  \label{tab:res_summary}
  \end{center}
 \vspace*{-.1in}
\end{table}

\paragraph{Nine Bar Truss} As shown in Fig \ref{fig:nineBar}, consider a truss with nine bars \citep{Mela07}, when the goal is to minimize two conflicting objectives: (a) the material volume of the truss ($J_1$), and (b) nodal displacement of the node $N$ ($J_2$):
\begin{align*}
&\min\left[\begin{array}{c}
J_{1}(\mathbf{x})\\
J_{2}(\mathbf{x})
\end{array}\right]=\left[\begin{array}{c}
L( \sum_{i=1}^{9} a_i x_i)\\
\frac{FL}{9E}\sum_{i=1}^{9} b_i / x_i)
\end{array}\right] \label{nineBar}  \tag{\textbf{E}$_2$} \\
& \text{subject to} \\
& x_i \in [l_i, 10],\qquad i=1,2,3 \\
& x_i \in \{1,5,10,15\},\qquad i=4,\ldots,9 \\
\end{align*}
The design variables, $\mathbf{x} = [x_1, \dots, x_9]$, are the cross-sectional areas of the bars, of which $x_1$, $x_2$, and $x_3$ are continuous and all others are discrete. The lower bounds for $x_1$ through $x_3$ are given as $l_1 = \frac{2}{3}$, $l_2 = \frac{1}{3}$, $l_3 = \frac{1}{3}$. $L$ is the length shown in Fig. \ref{fig:nineBar}, and $E$ is the Young's modulus of the bars. The coefficients $a_i$ and $b_i$ are $i^{\text{th}}$ elements of the sets $\mathbf{A}$ = \{1, 1, 1, $\sqrt{2}$, 1, $\sqrt{2}$, 1, $\sqrt{2}$, 1\} and $\mathbf{B}$ = \{4, 1, 1, $8\sqrt{2}$, 4, $2\sqrt{2}$, 4, $2\sqrt{2}$, 0\}, respectively. This example has a total of $|\mathcal{K}| = 4^6 = 4,096$ subproblems. The same set of discrete realizations as the the true one is identified; 78\% of the subproblems are pruned by Phase A and then another 20\% by Phase B (Table \ref{tab:res_summary}).

\section{Conclusions} \label{sec:concl}
A pruning scheme has been presented to reduce the number of discrete realizations to be considered for Pareto front generation of mixed-discrete bi-objective optimization.  Significant improvement in computational efficiency by the scheme has been verified on numerical examples.
Future work  includes extension to generic \textit{multi}-objective cases.

\newpage
\iftrue
\begin{acknowledgements}
This work was funded in part by research grant from the National Research Foundation of Korea (2011-0013956) and in part by the KI Project via KAIST Institute for Design of Complex Systems.
\end{acknowledgements}
\fi

\input{smo_pareto_hongahnchoi.bbl}

\end{document}

%% file: smo_pareto_hongahnchoi.bbl
\newcommand{\noopsort}[1]{} \newcommand{\printfirst}[2]{#1}
  \newcommand{\singleletter}[1]{#1} \newcommand{\switchargs}[2]{#2#1}